\documentclass[12pt]{amsart}
\usepackage{amssymb,a4wide, xy}
\xyoption{all}

\title {More on five commutator identities}
\author[G. DONADZE and M. LADRA]{G. DONADZE $^1$, M. LADRA $^2$}

\makeatletter
\let\c@equation=\c@subsection
\def\@begintheorem#1#2[#3]{\item[\normalfont\hskip\labelsep\bfseries
\noindent\textbf{(#2)}\ #1\unskip\@ifnotempty{#3}{\textup{\
(#3)}}.]\ignorespaces} \makeatother

\def\subsection#1{\refstepcounter{subsection}
\medskip\noindent
{\textbf{(\thesubsection)\ #1\unskip. }}\ignorespaces}
\def\subsubsection#1{\refstepcounter{subsubsection}
\smallskip\noindent
{\textbf{(\thesubsubsection)}\ \textit{#1\unskip. }}\ignorespaces}
\def\theequation{\thesection.\arabic{equation}}
\def\thesubsection{\thesection.\arabic{subsection}}
\def\thesubsubsection{\thesection.\arabic{subsection}.\arabic{subsubsection}}

\newtheorem{theorem}[subsection]{Theorem}

\newtheorem{definition}[subsection]{Definition}

\newtheorem{example}[subsection]{Example}

\newtheorem{lemma}[subsection]{Lemma}

\numberwithin{equation}{section}
\def\al{\alpha}
\def\be{\beta}

\def\de{\delta}

\def\ep{\epsilon}

\def\Gm{\Gamma}
\def\gm{\gamma}

\def\lra{\longrightarrow}

\def\otm{\otimes}
\def\ol{\overline}

\def\wt{\widetilde}

\def\uns{\underset}
\def\ovs{\overset}

\def\Ker{\operatorname{Ker}}
\def\Im{\operatorname{Im}}

\def\mod{\operatorname{mod}}

\def\lim{\operatorname{lim}}

\begin{document}
\maketitle \noindent {$^{1}$ \footnotesize \it A. Razmadze
Mathematical
Institute, Georgian Academy of Sciences, M. Alexidze St. 1,\\
0193 Tbilisi,
Georgia, e-mail: donad@rmi.acnet.ge}\\
{$^2$ \footnotesize \it Departamento de \'Algebra, Facultad de
Matem\'aticas, Universidad de Santiago de Compostela,\\ 15782
Santiago de Compostela, Spain, e-mail: ladra@usc.es}

\begin{abstract}
We prove that five well-known identities universally satisfied by
commutators in a group generate all universal commutator identities
for commutators of weight $4$.
\end{abstract}

\


\section*{Introduction}

\ For elements $x,y$ of a group we write $^xy=xyx^{-1}$ and $[x,y]=
xyx^{-1}y^{-1}$. The following commutator identities are universal
in the sense that they hold for any elements $x,y,z$ of an arbitrary
group:
$$
[x,x]=1,
$$
$$
[x,yz]=[x,y]\,^y[x,z] ,
$$
$$
[xy,z]=\, ^x[y,z][x,z],
$$
$$
[[y,x],^xz][[x,z],^zy][[z,y],^yx]=1,
$$
$$
^z[x,y]=[^zx,^zy].
$$
\noindent In \cite{El}  Ellis conjectured that, for any $n$, these
universal relations applied to commutators of weight $n$ generate
all universal relations between commutators of weight $n$. This
conjecture is stronger than Miller's result \cite{Ml}, who proved
that any universal relation among commutators is deduced from four
given ones without considering weights. Ellis considers his
conjecture as a nonabelian version of the Magnus-Witt theorem (see
\cite{M} and \cite{W}). To make his conjecture precise Ellis
introduced the structure of ``multiplicative Lie algebra''. Then
using the methods of homological algebra, he proved his conjecture
for $n=2$ and $n=3$.

This paper proves Ellis' conjecture for n=4 using essentially the
same tools.

\

 \ \


\section{Multiplicative Lie algebras}

\

This section is devoted to the formulation of Ellis' conjecture,
which he calls a nonabelian version of the Magnus-Witt theorem. We
first recall the notion of a multiplicative Lie algebra due to Ellis
\cite{El}.


\begin{definition}
A multiplicative Lie algebra consists of a multiplicative (possibly
nonabelian) group $L$ together with a binary function $\{ \;,\;
\}:L\times L\to L$, which we shall call Lie product, satisfying the
following identities for all $x,x',y,y',z$ in $L$

\begin{eqnarray}
\label{kuk2}\{x,x\} = 1\;, \\
\label{kuk3}\{x, yy'\} = \{x,y\}\;^y {\{x,y'\}}\;, \\
\label{kuk4}\{xx', y\} = \, ^x{\{x', y\}} \{x,y\}\;,\\
\label{kuk5}\{\{y,x\}, ^xz\}  \{\{x,z\}, ^zy\}  \{\{z,y\}, ^yx\} =
1\;,\\
\label{kuk6} ^z{\{x, y\}} = \{^z x,^z y\}\;.
\end{eqnarray}
\end{definition}

\medskip

In \cite{El} the following identities are deduced from
(\ref{kuk2})-(\ref{kuk6}):
\begin{eqnarray}
\label{kuk7}\{1,x\} = \{x, 1\} = 1\;, \\
\label{kuk8}\{y, x\} = \{x,y\}^{-1}\;, \\
\label{kuk9}^{\{x, y\}}{\{x', y'\}}= \, ^{[x, y]}{\{x', y'\}}\;,\\
\label{kuk10}\{[x,y], x'\} = [\{x,y\}, x']\;,\\
\label{kuk11}\{x^{-1}, y\} = {}^{x^{-1}}{\{x,y\}^{-1}}\quad
\text{and}\quad \{x, y^{-1}\} = {}^{y^{-1}}{\{x,y\}^{-1}}
\end{eqnarray}

\noindent for all $x, x', y, y'\in L$. Important examples of
multiplicative Lie algebras required for us are

\medskip


\begin{example}\label{kuk12}
Any group $P$ is a multiplicative Lie algebra with
$\{x,y\}=xyx^{-1}y^{-1}$ for all $x,y\in P$.
\end{example}

\medskip


\begin{example}
For any group $P$ there exists the {\it free} multiplicative Lie
algebra $\mathcal{L}(P)$ on $P$ which is characterized (up to
isomorphism) by the following two properties: $P$ is a subgroup of
$\mathcal{L}(P)$; and any group homomorphism $P\rightarrow L$ from
$P$ to a multiplicative Lie algebra $L$ extends uniquely to a
morphism of multiplicative Lie algebras $\mathcal{L}(P)\rightarrow
L$.
\end{example}

\medskip

The free multiplicative Lie algebra functor $\mathcal{L}$ is the
left adjoint of the forgetful functor from Multiplicative Lie
Algebras to Groups. The construction of $\mathcal{L}$ is given in
\cite{El} and more precisely in \cite{Bd}.

Let P be a group and $\Gm_n(P)$ be the subgroup of $\mathcal{L}(P)$
generated by the elements $\{\{\ldots\{\{x_1,x_2\},x_3\},\ldots
\},x_n\}$ for $x_i\in P$. In particular $\Gm_1(P)=P$. Then the group
identity morphism on $P$ induces a surjective morphism of
multiplicative Lie algebras
$$
\theta :\mathcal{L}(P)\twoheadrightarrow P
$$
in which $P$ has the structure of (\ref{kuk12}), and which restricts
to surjective group homomorphisms
$$
\theta_n :\Gm_n(P)\twoheadrightarrow \gm_n(P)
$$
for all $n\geq 1$, where $\gm_1(P)=P,\;\gm_n(P)=[\gm_{n-1}(P),P]$ is
the lower central series of $P$. Now we can exactly formulate the
Ellis' conjecture.

\medskip

{\bf Conjecture}. {\em If $P$ is a free group, then $\theta_n$ are
isomorphisms for all $n\geq 1$}.

\medskip

As we had already mentioned, the above  conjecture was proved in
\cite{El} for $n=2$ and $3$. The next section is devoted to the
proof for $n=4$.


\section{Ellis conjecture for commutators of weight 4}
\ We begin by recalling the notion of the nonabelian tensor product
introduced by Brown and Loday \cite{BL}  for a pair of groups
$G,\;H$ which act on themselves by conjugation and each of which
acts on the other compatibility, i.e.,
$$
^{(^{g}h)} g^{\prime}=\, ^{ghg^{-1}}g^{\prime},\quad
^{({^h}g)}h^{\prime}=\, ^{hgh^{-1}}h^{\prime}
$$
where $g,g^{\prime}\in G,\;h,h^{\prime}\in H$, and
$ghg^{-1},\;hgh^{-1}$ are elements of the free product $G*H$. The
nonabelian tensor product $G\otm H$ is the group generated by the
symbols $g\otm h$ subject to the relations
$$
gg^{\prime}\otm h=(^gg^{\prime}\otm \,^gh)(g\otm h)
$$
$$
g\otm hh^{\prime}=(g\otm h)(^hg\otm \,^hh^{\prime})
$$
\noindent for all $g,g^{\prime}\in G$ and $h,h^{\prime}\in H$.

 We will use the additive notations each time $G\otimes H$ is abelian.

In the sequel, unless specified, the tensor product of groups
$G\otimes H$ belongs to three kinds for which the compatibility
conditions hold:

\noindent(1) $G$ is a normal subgroup of $H$ and actions are given
by conjugations;

\noindent(2) $G$ is an abelian quotient of some normal subgroup of
$H$, the action of $H$ on $G$ is induced by conjugation and the
action of $G$ on $H$ is trivial;

\noindent(3) $H=P_{ab}$ and $G$ is a  quotient of $[P,P]_{ab}$, for
some group $P$, the action of $H$ on $G$ is induced by conjugation
and the action of $G$ on $H$ is trivial.

\medskip

Let $P$ be a group. Define $[P,P]_{ab}\otimes P_{ab}$ according to
(3).  As $[P,P]_{ab}$ is a $P_{ab}$-module, \cite[Proposition
3.2]{Gu} says that $[P,P]_{ab}\otimes P_{ab}$ is isomorphic to
$[P,P]_{ab}\otimes_{P_{ab}}IP_{ab}$, where $IP_{ab}$ denotes the
augmentation ideal of $P_{ab}$. Hence $[P,P]_{ab}\otimes P_{ab}$ is
abelian.

\medskip

\begin{lemma}\label{1}
Let $P$ be a group. Then we have the following equalities  in
$[P,P]_{ab}\otimes P_{ab}$:
\begin{eqnarray}
&([x,y][x',y'])\otimes z = [x,y]\otimes z + [x',y']\otimes z\;, \\
\label{kuk2.3}&[[a,b],y]\otimes x + [x,[a,b]]\otimes y = 0\;, \\
&[{}^p z,x]\otimes y + [y,{}^p z]\otimes x - [z,x]\otimes y -
[y,z]\otimes x = 0\;.
\end{eqnarray}
for any $a,b,x,y,z,p \in P$.
\end{lemma}
\begin{proof} We only prove the second and third equalities. In fact,
\begin{align*}
&[[a,b],y]\otimes x +[x,[a,b]]\otimes y=[a,b] \,
{}^{y}{[b,a]}\otimes x + \, {}^{x}{[a,b]}[b,a]\otimes y \\   & =
[a,b]\otimes x    +  \ {}^{y}{[b,a]}\otimes x+ \, ^x[a,b]\otimes
y+[b,a]\otimes y \\  & =[a,b]\otimes x + \, {}^{x}{[a,b]}\otimes {}^
xy + [b,a]\otimes y
 + \ {}^{y}{[b,a]}\otimes {}^y x \\  & = [a,b]\otimes xy +[b,a]\otimes
yx = [a,b]\otimes xy + [b,a]\otimes xy = 0\;,
\end{align*}
and
\begin{align*}
& [{}^p z,x]\otimes y + [y,{}^p z]\otimes x - [z,x]\otimes y -
[y,z]\otimes x \\  & = [[p,z]z,x]\otimes y + [y,[p,z]z]\otimes x   -
\ [z,x]\otimes y - [y,z]\otimes x \\  & =([z,x][[p,z],x])\otimes y +
([y,[p,z]][y,z])\otimes x - [z,x]\otimes y    - \ [y,z]\otimes x \\
& = [z,x]\otimes y + [[p,z],x]\otimes y + [y,[p,z]]\otimes x +
[y,z]\otimes x - [z,x]\otimes y    - \ [y,z]\otimes x \\ &
=[[p,z],x]\otimes y + [y,[p,z]]\otimes x = 0 \;, \text{by}\;
(\ref{kuk2.3}).
\end{align*}
\end{proof}
\medskip
Let $P$ be a group. $P$  and $P\otimes P$ are $P$-crossed modules
and they act on each other via their images in the basis $P$, i.e.,
$$
^z(x\otimes y)= \, ^zx\otimes {}^zy, \quad ^{x\otimes y}z= \,
^{[x,y]}z \;.
$$
Thus, in the next lemma, $(P\otimes P,P)$ is a pair equipped with
compatible actions and we can define the nonabelian tensor product
$(P\otimes P)\otimes P$. In order to describe $(P\otimes P)\otimes
P$ more precisely, assume that $F$ is the free group generated by
symbols $x\otimes y$ , $x,y\in P$. Then,  $(P\otimes P)\otimes P$
will be the group generated by symbols $f\otimes z$, $f\in F$, $z\in
P$, subject to the following relations
\begin{eqnarray*}
&(ff'\otimes z)(f\otimes z)^{-1}(^ff'\otimes {}^f z)^{-1} = 1,\\
&(f\otimes zz')({}^z f\otimes {}^{z} {z'})^{-1}(f\otimes
z)^{-1} = 1,\\
&\overline{f}\otimes z = 1, \overline{f}\in \overline{F},
\end{eqnarray*}
where $f,f'\in F$, $z,z'\in P$, $f$ acts on $z$ via its image in
$P\otm P$ and $\overline{F}$ is the normal subgroup of $F$ generated
by the following elements
\begin{eqnarray*}
& (xx'\otimes y)(x\otimes y)^{-1}({}^{x}{x'}\otimes {}^x
y)^{-1}\;,\\
& (x\otimes y y')({}^y x\otimes {}^{y}{y'})^{-1}(x\otimes y)^{-1}\;.
\end{eqnarray*}

\medskip

\begin{lemma}\label{2}
Assume that $P$ is a group, $F$ is the aforementioned group, i.e.,
$F$ is the free group generated by symbols $x\otm y$, $x,y\in P$ and
 $[P,P]_{ab} \otimes P_{ab}$ is defined
as in Lemma\;(\ref{1}). Then there is a well-defined homomorphism
$\de :(P\otimes P)\otimes P \lra [P,P]_{ab} \otimes P_{ab}$, given
as follows: if $f=\uns{i}{\prod}(x_i\otimes y_i)^{\ep _i}\in F$,
$\ep_i=\pm1$, and $z\in P$ then
$$
f\otimes z\mapsto \uns{i}{\sum}\ep _i([x_i,y_i]\otimes z
+[z,x_i]\otimes y_i +[y_i,z]\otimes x_i)\;,
$$
where $f$ is identified with its image into $P\otm P$.
\end{lemma}
\begin{proof} Taking into account the relations above, we have to check
the following equalities:
\begin{align}\label{d1}
\de ((ff^{\prime}\otimes z)(f\otimes z)^{-1}(^ff^{\prime}\otimes
{}^fz)^{-1})=0,
\end{align}
\begin{align}\label{d2}
\de ((f\otimes zz^{\prime})(^zf\otimes {}^zz^{\prime})^{-1}(f\otimes
z)^{-1})=0,
\end{align}
\begin{align}\label{d4}
\de (((xx^{\prime}\otimes y)(x\otimes y)^{-1}(^xx^{\prime}\otimes {}
^xy)^{-1})\otimes z)=0,
\end{align}
\begin{align}\label{d5}
\de (((x\otimes yy^{\prime})(^yx\otimes {}
^yy^{\prime})^{-1}(x\otimes y)^{-1})\otimes z)=0,
\end{align}
\

\noindent where $f,f^{\prime}\in F$ and
$x,y,x^{\prime},y^{\prime},z\in P$.

 The proof of (\ref{d1}) will be
trivial, if we show that $\de (^{f}f'\otimes \;^{f}z)=\de (f'\otimes
z)$ for all $f,f'\in F$ and $z\in P$. It suffices to take $f=x\otm
y$ and $f'=x'\otm y'$, where $x,y,x',y'\in P$. Thus, we need to show
that $\de (^{(x\otm y)}(x'\otm y')\otm ^{(x\otm y)}z)=\de ((x'\otm
y')\otm z)$, which is equivalent to the equality $\de ((x'\otm
y')\otm ^{[x, y]}z)=\de ((x'\otm y')\otm z)$. Clearly $\de ((x'\otm
y')\otm z)=\de ((^{[x,y]}x^{\prime}\otimes
{}^{[x,y]}y^{\prime})\otimes \,^{[x,y]}z)$, hence we have to check
that $\de ((x'\otm y')\otm ^{[x,y]}z)=\de
((^{[x,y]}x^{\prime}\otimes {}^{[x,y]}y^{\prime})\otimes
\,^{[x,y]}z)$, which is equivalent to the following:
$$
\de (((x^{\prime}\otimes y^{\prime})(^{[x,y]}x^{\prime}\otimes {}
^{[x,y]}y^{\prime})^{-1})\otimes z)=0\;.
$$
One has:
\begin{align*}
& \de (((x^{\prime}\otimes y^{\prime})(^{[x,y]}x^{\prime}\otimes {}
^{[x,y]}y^{\prime})^{-1})\otimes z)
\\ & = [x^{\prime},y^{\prime}]\otimes z+[z,x^{\prime}]\otimes
y^{\prime}+[y^{\prime},z]\otimes x^{\prime} -[{} ^{[x,y]}x^{\prime},
{} ^{[x,y]}y^{\prime}]\otimes z  - [z, {} ^{[x,y]}x^{\prime}]\otimes
{} ^{[x,y]}y^{\prime} \\ &  \qquad  \qquad   - [{}
^{[x,y]}y^{\prime},z]\otimes {} ^{[x,y]}x^{\prime} \\ & = [{}
^{[x,y]}x^{\prime},{} ^{[x,y]}y^{\prime}]\otimes {} ^{[x,y]}z+
[{}^{[x,y]}z, {} ^{[x,y]}x^{\prime}]\otimes {} ^{[x,y]}y^{\prime}+
[{}^{[x,y]}y^{\prime}, {} ^{[x,y]}z]\otimes {}^{[x,y]}x^{\prime}
\\ & \qquad   \qquad -[{} ^{[x,y]}x^{\prime},{} ^{[x,y]}y^{\prime}]\otimes z- [z,{}
^{[x,y]}x^{\prime}]\otimes {} ^{[x,y]}y^{\prime}  - [{}
^{[x,y]}y^{\prime},z]\otimes {} ^{[x,y]}x^{\prime}\\ & =
[{}^{[x,y]}z,{} ^{[x,y]}x^{\prime}]\otimes {} ^{[x,y]}y^{\prime}+
[{}^{[x,y]}y^{\prime},{} ^{[x,y]}z]\otimes {} ^{[x,y]}x^{\prime}-
[z,{}^{[x,y]}x^{\prime}]\otimes {} ^{[x,y]}y^{\prime}- [{}
^{[x,y]}y^{\prime},z]\otimes {} ^{[x,y]}x^{\prime} \\ & = 0 \, , \,
\text{by previous lemma}.
\end{align*}

 We will check only
(\ref{d2}) and (\ref{d4}), because (\ref{d5}) is similar to
(\ref{d4}).

(\ref{d2}): One can easily see that it is enough to consider
$f=x\otimes y$.
\begin{align*}
& \de (((x\otimes y)\otimes zz^{\prime})((^zx\otimes {}^zy)\otimes
{} ^zz^{\prime})^{-1}((x\otimes y)\otimes z)^{-1})\\ & =
[x,y]\otimes zz^{\prime}+[zz^{\prime},x]\otimes
y+[y,zz^{\prime}]\otimes x-[^zx,^zy]\otimes {} ^zz^{\prime}-
[^zz^{\prime},^zx]\otimes {}^zy \\ & \qquad
-[^zy,^zz^{\prime}]\otimes {} ^zx - \ [x,y]\otimes z-[z,x]\otimes y
-[y,z]\otimes x \\ & =[x,y]\otimes z + [^zx,^zy]\otimes {}
^zz^{\prime}+[^zz^{\prime},^zx]\otimes y+ [z,x]\otimes y + \
[y,z]\otimes x+[^zy,^zz^{\prime}]\otimes x  \\ &  \qquad -
[^zx,^zy]\otimes {} ^zz^{\prime}- [^zz^{\prime},^zx]\otimes {}^zy-
[^zy,^zz^{\prime}]\otimes {} ^zx-[x,y]\otimes z -[z,x]\otimes
y-[y,z]\otimes x \\ & =0 \, , \; \text{since} \ [a,b]\otimes {}
^zp=[a,b]\otimes p.
\end{align*}

(\ref{d4}):
\begin{align*}
& \de (((xx^{\prime}\otimes y)(x\otimes y)^{-1}(^xx^{\prime}\otimes
{} ^xy)^{-1})\otimes z) \\ & = [xx^{\prime},y]\otimes
z+[z,xx^{\prime}]\otimes y +[y,z]\otimes x x^{\prime}-[x,y]\otimes
z-[z,x]\otimes y-[y,z]\otimes x \\ &  \qquad - \
[^xx^{\prime},^xy]\otimes
z  -[z,^xx^{\prime}]\otimes {} ^xy- [^xy,z]\otimes {} ^xx^{\prime}\\
& =[^xx^{\prime},^xy]\otimes z+ [x,y]\otimes z+[z,x]\otimes y + \
[^xz,^xx^{\prime}]\otimes y+[y,z]\otimes x  +[^xy,^xz]\otimes {}
^xx^{\prime} \\
& \qquad -[x,y]\otimes z-[z,x]\otimes y-[y,z]\otimes x - \
[^xx^{\prime},^xy]\otimes
z-[z,^xx^{\prime}]\otimes {} ^xy- [^xy,z]\otimes {} ^xx^{\prime}\\
& =[^xz,^xx^{\prime}]\otimes {} ^xy+ [^xy,^xz]\otimes {}
^xx^{\prime} - \ [z,^xx^{\prime}]\otimes {}^xy- [^xy,z]\otimes {}
^xx^{\prime} \\
& = 0 \,, \, \text{by the previous lemma}.
\end{align*}
\end{proof}

\medskip

\begin{lemma}\label{3}
Let $P$ be a group and $H_2(P)=0$. Then the homomorphism $\de
:(P\otimes P)\otimes P\lra [P,P]_{ab}\otimes P_{ab}$\; introduced in
Lemma (\ref{2}) factors through
$(\gamma_2(P)/\gamma_3(P))\otimes_{\mathbb{Z}}P_{ab}$. Thus, $\de ^*
:(\gamma_2(P)/\gamma_3(P))\otimes_{\mathbb{Z}}P_{ab} \lra
[P,P]_{ab}\otimes P_{ab}$ given by
$$
[x,y]\otimes z \mapsto [x,y]\otimes z+[z,x]\otimes y+[y,z]\otimes x
$$
is well defined.
\end{lemma}
\begin{proof} By \cite{BL} we have $[P,P] \cong (P\otimes P)/X_1$, where $X_1$
is the normal subgroup $P\otimes P$ generated  by $x\otimes x$, for
all $x\in P$. Therefore,
$$
(\gm _2(P)/\gm _3(P)) \otimes_\mathbb{Z} P_{ab} \cong ((P\otimes
P)\otimes P)/X_2
$$
where $X_2$  is the normal subgroup of $(P\otimes P)\otimes P$
generated by $(x\otimes x)\otimes z$, $([x,y]\otimes z)\otimes
p,\;(x\otimes y)\otimes [p,q]$ for all $x,y,z,p,q,\in P$. By the
previous lemma it is enough to check the following:
$$
\de ((x\otimes x)\otimes z)=0,
$$
$$
\de (([x,y]\otimes z)\otimes p))=0,
$$
$$
\de ((x\otimes y)\otimes [p,q])=0.
$$
We have:
$$
\de ((x\otimes x)\otimes z)=[x,x]\otimes z+[z,x]\otimes
x+[x,z]\otimes x=0,
$$
$$
\de (([x,y]\otimes z)\otimes p))=[[x,y],z]\otimes p+[p,[x,y]]\otimes
z+[z,p]\otimes [x,y]=0 \;, \text{by Lemma}\;(\ref{1}),
$$
$$
\de ((x\otimes y)\otimes [p,q])=[x,y]\otimes [p,q]+[[p,q],x]\otimes
y+ [y,[p,q]]\otimes x=0 \;, \text{by Lemma}\;(\ref{1}).
$$
\end{proof}

\medskip

Given a group $P$, let $\al :[P,P]\otimes P \lra \gm _3(P)$ be the
commutator map: $\al ([x,y]\otimes z)=[[x,y],z]$. $\al$ induces the
following homomorphisms:
$$
\al_1 :[P,P]_{ab}\otimes P_{ab}\lra \gm _3(P)/[[P,P],[P,P]],$$
$$
\al_2 :(\gm_2(P)/\gm _3(P))\otimes P_{ab} \lra \gm _3(P)/\gm _4(P).
$$
Define $i:\Ker\al \rightarrow \Ker\al_1$ and $i':\Ker\al_1
\rightarrow \Ker\al_2$ as restrictions of the natural projections
$[P,P]\otimes P \lra [P,P]_{ab}\otimes P_{ab}$ and
$[P,P]_{ab}\otimes P_{ab}\lra (\gm_2(P)/\gm _3(P))\otimes P_{ab}$,
respectively.

\medskip

\begin{lemma}\label{4}
Let $P$ be a free group and  define $\de^*$  as in Lemma (\ref{3}).
Then $\Ker \al_1=\Im \de^*$ and  $i':\Ker\al_1 \lra \Ker\al_2$ is an
isomorphism.
\end{lemma}
\begin{proof} Part one: Using diagram chasing we easily check that
$i:\Ker\al \lra \Ker\al_1$ is surjective. \cite[Theorem 9]{El} says
that $\Ker \al$ is generated by $([x,y]\otimes
\;^yz)([y,z]\otimes\;^zx)([z,x]\otimes \;^xy) \in [P,P]\otimes P$
and $p\otimes p \in [P,P]\otimes P$ for all $x,y,z\in P$ and $p\in
[P,P]$. Since $i$ is surjective, the generators of $\Ker\al_1$ will
be  $[x,y]\otimes z+[z,x]\otimes y+[y,z]\otimes x \in
[P,P]_{ab}\otimes P_{ab}$ for all $x,y,z\in P$. Hence $\Ker
\al_1=\Im \de^*.$

Part two: Using diagram chasing we easily check that $i':\Ker\al_1
\lra \Ker\al_2$ is surjective. To show the injectivity we need to
check the following: if $\omega \in \Ker\al_1$, then $\de
^*i'(\omega )= 3\omega $. In fact, by discussion above it suffices
to take $w=[x,y]\otm z+[z,x]\otm y+[y,z]\otm x$. We have
\begin{align*}
& \de^*i'(w)=\de^*([x,y]\otm z)+\de^*([z,x]\otm y)+\de^*([y,z]\otm
x) \\& = [x,y]\otm z+[z,x]\otm y+[y,z]\otm x + [z,x]\otm y+[y,z]\otm
x \\& \qquad  \qquad  +[x,y]\otm x + [y,z]\otm x+[x,y]\otm
z+[z,x]\otm y \\& = 3([x,y]\otm z+[z,x]\otm y+[y,z]\otm x)=3w.
 \end{align*}
 Therefore,
for injectivity of $i'$, it is sufficient to show that $\Ker\al_1$
is torsion free. We have $\Ker\al_1 \cong H_3(P_{ab})$ (see
\cite[Theorem 6.7]{IH1} and \cite{El}). Since $P$ is free,
$H_3(P_{ab})\cong P_{ab}\wedge P_{ab}\wedge P_{ab}$ (see
\cite{Bro,G,HS}) which is torsion free.
\end{proof}

\medskip

Given a group $P$, we have the short exact sequence of groups
\begin{align}\label{m1}
1\lra \gm_3(P)/[\gm_2(P),\gm_2(P)]\lra [P,P]_{ab}\lra
\gm_2(P)/\gm_3(P)\lra 1.
\end{align}
Suppose $X$ be one of the groups in the sequence (\ref{m1}). Define
$X\otimes P$ according to (2). Assume that the homomorphisms
$$
\be_0:(\gm_3(P)/[\gm_2(P),\gm_2(P)])\otm P\lra
\gm_4(P)/[\gm_3(P),\gm_2(P)]$$
$$
\be_1:[P,P]_{ab}\otm P\lra \gm_3(P)/[\gm_3(P),\gm_2(P)],$$
$$
\be_2:(\gm_2(P)/\gm_3(P))\otm P\lra \gm_3(P)/\gm_4(P),
$$
are defined by taking commutators. These homomorphisms are well
defined because the restrictions to $\gamma_2(P)$ of the actions of
$P$ on $\gm_3(P)/[\gm_3(P),\gm_2(P)]$,
$\gm_4(P)/[\gm_3(P),\gm_2(P)]$ and $\gamma_3(P)/\gamma_4(P)$ induced
by conjugation are trivial. If $P$ is a free group, then there is a
short exact sequence of groups
\begin{align}\label{m2}
0\lra \Ker\be_0 \lra \Ker\be_1 \lra \Ker\be_2 \lra 0.
\end{align}
In fact, since the groups in Sequence (\ref{m1}) are $P$-modules and
act trivially on $P$, we have $\otm P \cong \otm _PIP$, where
$IP=\Ker (\mathbb{Z}(P)\rightarrow \mathbb{Z})$. Since $P$ is free,
$IP$ is a free $P$-module. Therefore, $(\ref{m1})\otm P$ is a short
exact sequence. This implies that (\ref{m2}) is a short exact
sequence.

\medskip

\begin{lemma}\label{5} Let P be a free group. Assume that

\noindent 1) $A$ is a normal subgroup of
$(\gm_3(P)/[\gm_2(P),\gm_2(P)])\otm P$ generated by all $x\otm y \in
(\gm_3(P)/[\gm_2(P),\gm_2(P)])\otm P$, where $y\in [P,P]$;

 \noindent 2) $A'$ is the natural image of $A$ into $[P,P]_{ab}\otm P$;

\noindent 3) $B$ is a normal subgroup of $[P,P]_{ab}\otm P$
generated by all $z\otm z \in [P,P]_{ab}\otm P$ where $z\in [P,P]$
and in $z\otm z$, the first $z$ is identified with its image in
$[P,P]_{ab}$.

 \noindent Then

 \noindent $\mathbf{(a)}$
$\Ker\be_2 \cong \Ker\be_1/(A'+B)$.

 \noindent $\mathbf{(b)}$ $\Ker\be_0$
is generated by $A$ and the set of elements \;\;$[y,{}
^{y^{-1}}[x,y]]\otm x+[{} ^{y^{-1}}[x,y],x]\otm {} ^xy$, \;
$[y,{}^{y^{-1}}[p,q]]\otm x+[{} ^{y^{-1}}[p,q],x]\otm {} ^xy+ [q,{}
^{q^{-1}}[x,y]]\otm p+[{} ^{q^{-1}}[x,y],p]\otm {}^pq$ \;\;for all
$x,y,p,q\in P$.
\end{lemma}
\begin{proof}$\mathbf{(a)}$:
Denote by $j$ the natural homomorphism $\Ker\be_1\lra \Ker\be_2$.
Since $\gm_2(P)/\gm_3(P)$ and $P$ act trivially on each other, by
\cite{BL} we have
\begin{align}\label{m3}
(\gm_2(P)/\gm_3(P))\otm P \cong (\gm_2(P)/\gm_3(P))\otm P_{ab}\;.
\end{align}
Thanks to this we easily see that $j$ sends $(A'+B)$ to zero and
induces a homomorphism $j^*:\Ker\be_1/(A'+B)\lra \Ker\be_2$. Assume
that $[P,P]\wedge P$ is defined naturally (i.e. as in \cite{El}) and
$\al ^*: [P,P]\wedge P\lra \gm_3(P)$ is the commutator map. We have
the natural projection
$$
[P,P]\wedge P\lra ([P,P]_{ab}\otm P)/ (A'+B),
$$
which induces a homomorphism:
$$
\Ker\al^*\ovs{\tau}{\lra} \ \Ker\be_1/(A'+B).
$$
Using diagram chasing we easily check that  $\tau $ is an
epimorphism. Hence the composition $ \Ker\al^*
\ovs{\tau}{\lra}\Ker\be_1/(A'+B)\ovs{j^*}{\lra} \Ker\be_2 $ is an
epimorphism. Prove that $j^*\circ\tau $ is an isomorphism. Since $P$
is a direct limit of its finitely generated subgroups and this
system is compatible with $j^*, \tau, \al^*$ and $\beta_2$, without
lost of generality we can assume that $P$ is a free group with
finite basis. Then $H_3(P_{ab})\cong P_{ab}\wedge P_{ab}\wedge
P_{ab}\cong \mathbb{Z}^n$ for some $n\in \mathbb{N}$. Moreover,
taking into account (\ref{m3}) and the previous lemma, we have
$$
\Ker\be_2\ \cong \Ker\al_2 \cong \Ker\al_1 \cong H_3(P_{ab})\cong
\mathbb{Z}^n,
$$
where $\al_1$ and $\al_2$ are defined as above. On the other hand
there is an isomorphism $\Ker\al^*\cong H_3(P_{ab})$ (see
\cite{El}). Thus, both of $\Ker\al^*$ and $\Ker\beta_2$ are
isomorphic to $\mathbb{Z}^n$. Therefore, any epimorphism $\Ker\al^*
\to \Ker\beta_2$ (in particular $j^*\circ\tau $) is an isomorphism.
Hence $j^*$ is an isomorphism and we have $\mathbf{(a)}$ .

$\mathbf{(b)}$: Note that $B$ is generated by all $[x,y]\otm
[x,y]\in [P,P]_{ab}\otm P$ and $([p,q]\otm [x,y]+[x,y]\otm [p,q])\in
[P,P]_{ab}\otm P $, where $x,y,p,q\in P$. Therefore, taking into
account (\ref{m2}) and $\mathbf{(a)}$, it is sufficient to prove
that in $[P,P]_{ab}\otm P$ the following hold:
$$
[x,y]\otm [x,y]=[y,{} ^{y^{-1}}[x,y]]\otm x+[^{y^{-1}}[x,y],x]\otm
{} ^xy,
$$
\begin{align*}
&  \qquad \qquad [p,q]\otm [x,y]+[x,y]\otm [p,q] \\ & =[y,{}
^{y^{-1}}[p,q]]\otm x+[^{y^{-1}}[p,q],x] \otm {} ^xy + \ [q,{}
^{q^{-1}}[x,y]]\otm p+[^{q^{-1}}[x,y],p]\otm {} ^pq,
\end{align*}
 for any $x,y,p,q\in P$. Both of these equalities will be
clear, if we prove the following:
$$
[p,q]\otm [x,y]=[y,{} ^{y^{-1}}[p,q]]\otm x+ [^{y^{-1}}[p,q],x]\otm
{} ^xy.
$$
for all $x,y,p,q\in P$. We have:
\begin{align*}
& [p,q]\otm [x,y]=[p,q]\otm x \, ^yx^{-1}=[p,q]\otm x+ {}
^x[p,q]\otm {}^{xy}x^{-1} \\ & =[p,q]\otm x + \
(^{xyy^{-1}}[p,q]\otm {}
^{xy}x^{-1}+ \, ^{y^{-1}}[p,q]\otm xy)- \, ^{y^{-1}}[p,q]\otm xy \\
& =[p,q]\otm x + \ ^{y^{-1}}[p,q]\otm (xy)x^{-1}- \,
^{y^{-1}}[p,q]\otm xy \\ & =[p,q]\otm x+ \, ^{y^{-1}}[p,q]\otm {}
^xy - \ ^{y^{-1}}[p,q]\otm x- \, ^{xy^{-1}}[p,q]\otm {} ^xy \\ & =
([p,q]\otm x-\, ^{y^{-1}}[p,q]\otm
x) + \ (^{y^{-1}}[p,q]\otm {}^xy- \, ^{xy^{-1}}[p,q]\otm {}^xy) \\
& = (^{yy^{-1}}[p,q]\otm x+ \, ^{y^{-1}}[p,q]^{-1}\otm x) + \
(^{y^{-1}}[p,q]\otm {} ^xy+ \, ^{xy^{-1}}[p,q]^{-1}\otm {} ^xy)
\\ & = [y,{} ^{y^{-1}}[p,q]]\otm x+[^{y^{-1}}[p,q],x]\otm {} ^xy.
\end{align*}
\end{proof}

\medskip

\begin{lemma}\label{6}
Let $P$ be a free group and let $\al ^{\prime}:\gm_3(P)\otm P\lra
\gm_4(P)/[\gm_3(P),\gm_2(P)]$ be the homomorphism defined by taking
commutators, i.e., $[[x,y],z]\otm p\mapsto [[[x,y],z],p]$, for
$x,y,z,p\in P$. Then $\Ker\al^{\prime}$ is generated by the
subgroups $[\gm_2(P),\gm_2(P)]\otm P$ and $\gm_3(P)\otm \gm_2(P)$
and the set of elements $([y, {} ^{y^{-1}}[x,y]]\otm
x)([^{y^{-1}}[x,y],x]\otm \, ^xy)$ and $([y,{}^{y^{-1}}[p,q]]\otm
x)([{}^{y^{-1}}[p,q],x]\otm \, ^xy) ([q,{}^{q^{-1}}[x,y]]\otm
p)([^{q^{-1}}[x,y],p]\otm {} ^pq)$, for all $x,y,p,q\in P$.
\end{lemma}
\begin{proof} Assume that $A$ is defined as in Lemma (\ref{5}).
Then we have an isomorphism
$$
\frac{\gamma_3(P)\otm P}{\gamma_3(P)\otm \gamma_2(P)+
[\gamma_2(P),\gamma_2(P)]\otm P}\cong
\frac{(\gamma_3(P)/[\gamma_2(P),\gamma_2(P)])\otm P}{A}
$$
given in a natural way. The rest of the proof is a consequence of
Lemma (\ref{5})$\mathbf{(b)}$.
\end{proof}

\medskip

\begin{theorem}\label{7}
If  $P$ is a free group, then $\theta_4:\Gm_4(P)\lra \gm_4(P)$ is an
isomorphism.
\end{theorem}
\begin{proof}Since surjectivity of $\theta_4$ is obvious, we will prove the
injectivity.

Let $\Gm_4^{\prime}(P)$ be the subgroup of $\Gm_4(P)$ generated by
$\{\{\{p_1,p_2\},[p,q]\},p_3\}$ and $\{\{\{p_1,p_2\},p_3\},[p,q]\}$,
for all $p_1,p_2,p_3,p,q\in P$. We easily see that
$\Gm_4^{\prime}(P)$ is a normal subgroup of $\mathcal{L}(P)$ and
$\theta_4(\Gm_4^{\prime}(P))=[\gm_3(P),\gm_2(P)]$. Define the
homomorphisms $\ol{\theta_4}$ and $\wt{\theta_4}$:
$$
\Gm_4^{\prime}(P)\xrightarrow{ \ \ol{\theta_4} \ } \theta_4
(\Gm_4^{\prime}(P))=[\gm_3(P),\gm_2(P)]\;\; \text{is the restriction
of}\;\;\theta_4 \;\;\text{on}\;\;\Gm_4^{\prime}(P)\;;$$
$$
\Gm_4(P)/\Gm_4^{\prime}(P)\xrightarrow{ \ \wt{\theta_4} \ }
\gm_4(P)/\theta_4(\Gm_4^{\prime}(P))=\gm_4(P)/[\gm_3(P),\gm_2(P)]\;\;\text{is
induced by}\;\;\theta_4\;.
$$
\noindent The proof of (\ref{7}) will be done, if we show that
$\ol{\theta_4}$ and $\wt{\theta_4}$ are injective. Using
(\ref{kuk9}) and (\ref{kuk10}) we easily show that
$\Gm_4^{\prime}(P)\subset \Gm_3(P)$. Since $\theta_3$ is an
isomorphism, $\ol{\theta_4}$ will be injective. In order to show
injectivity of $\wt{\theta_4}$, we construct the homomorphism
$$
 \hspace{-4cm} (\theta_3^{-1}\wt{\otm }P):\gm_3(P)\otm P\lra
\Gm_4(P)\;,$$
$$ \qquad \qquad \qquad  [[x,y],z]\otm p\mapsto
\{\theta_3^{-1}([[x,y],z]),p\}=\{\{\{x,y\},z\},p\}.
$$
Taking into account (\ref{kuk9}), it is trivial to check that
$\theta_3^{-1}\wt{\otm }P$ is well defined. Then, the following
composition
$$
\gm_3(P)\otm P\xrightarrow{ \ \theta_3^{-1}\wt{\otm \ } P}
\Gm_4(P)/\Gm_4^{\prime}(P)\ovs{\wt{\theta_4}}{\lra}
\gm_4(P)/[\gm_3(P),\gm_2(P)]
$$
is the map $\al^{\prime} :\gm_3(P)\otm P\lra
\gm_4(P)/[\gm_3(P),\gm_2(P)]$ defined in Lemma (\ref{6}). Since
$\theta_3^{-1}\wt{\otm }P$ is onto,
$\Ker\wt{\theta_4}=(\theta_3^{-1}\wt{\otm } P)(\Ker\al^{\prime} )$.
Hence, the generators of $\Ker\wt{\theta_4}$ are the images by
$\theta_3^{-1}\wt{\otm } P$ of the set of generators given in Lemma
(\ref{6}). Thus, we have to show the following:
\begin{align}\label{m4}
(\theta_3^{-1}\wt{\otm} P)([\gm_2(P),\gm_2(P)]\otm P)\subset
\Gm_4^{\prime}(P),
\end{align}
\begin{align}\label{m5}
(\theta_3^{-1}\wt{\otm } P)(\gm_3(P)\otm \gm_2(P))\subset
\Gm_4^{\prime}(P),
\end{align}
\begin{align}\label{m6}
(\theta_3^{-1}\wt{\otm } P)(([y,{} ^{y^{-1}}[x,y]]\otm
x)([^{y^{-1}}[x,y],x]\otm {} ^xy))\in \Gm_4^{\prime}(P),
\end{align}
\begin{align}\label{m7}
(\theta_3^{-1}\wt{\otm } P)(([y,{} ^{y^{-1}}[p,q]]\otm
x)([^{y^{-1}}[p,q],x]\otm {} ^xy) ([q,{} ^{q^{-1}}[x,y]]\otm
p)([^{q^{-1}}[x,y],p]\otm {} ^pq))\in \Gm_4^{\prime}(P).
\end{align}
(\ref{m4}) and (\ref{m5}) are trivial inclusions. For (\ref{m6}) and
(\ref{m7}), note that there are the following congruences $\mod \
\Gm_4^{\prime}(P)$:
$$
\{\{\{z_1,z_2\},^zz^{\prime}\},z_3\}\equiv
\{\{\{z_1,z_2\},z^{\prime}\},z_3\},
$$
$$
\{\{\{z_1,z_2\},z_3\},^zz^{\prime}\}\equiv
\{\{\{z_1,z_2\},z_3\},z^{\prime}\},
$$
$$
\{\{\{z_1,z_2\},z_3\},^{\{z,z^{\prime}\}}z_4\}\equiv
\{\{\{z_1,z_2\},z_3\},z_4\},
$$
for all $z_1,z_2,z_3,z_4,z,z^{\prime}\in P$. These relations and
(\ref{kuk5}) imply that
\begin{align*}
&(\theta_3^{-1}\wt{\otm } P)(([y,{} ^{y^{-1}}[x,y]]\otm
x)([^{y^{-1}}[x,y],x]\otm {} ^xy))\\ & = \{\{y,^{y^{-1}}\{x,y\}\},
x\}\{\{^{y^{-1}}\{x,y\},x\}, ^xy\}\equiv
(\{\{x,y\},\{x,y\}\})^{-1}=1,
\end{align*}
\begin{align*} & (\theta_3^{-1}\wt{\otm } P)(([y,{}
^{y^{-1}}[p,q]]\otm x)([^{y^{-1}}[p,q],x]\otm {} ^xy) ([q,{}
^{q^{-1}}[x,y]]\otm p)([^{q^{-1}}[x,y],p]\otm {}^pq))
\\ &
=\{\{y,^{y^{-1}}\{p,q\}\},x\}\{\{^{y^{-1}}\{p,q\},x\}, ^xy\}
\{\{q,^{q^{-1}}\{x,y\}\},p\}\{\{^{q^{-1}}\{x,y\},p\},^pq\}\}
\\ &
\equiv (\{\{x,y\},\{p,q\}\})^{-1} (\{\{p,q\},\{x,y\}\})^{-1}=1,
\end{align*}
 where the congruences being still taken $\mod \
\Gm_4^{\prime}(P)$. Thus, (\ref{m6}) and (\ref{m7}) are proved.
\end{proof}

\ \
\begin{center}
{\bf Acknowledgements}
\end{center}
The second author was supported by Xunta de Galicia,
PGIDIT06PXIB371128PR and the MEC (Spain), MTM 2006-15338-C02-01
(European FEDER support included). We are grateful to Professor Nick
Inassaridze for discussion during the period of development of the
paper. We are sorry that he refused to be one of the coauthors of
this paper. We also thank the referees for several helpful
suggestions which have significantly contributed to improve the
paper.
 \
 \

\end{document}